\documentclass[12pt]{article}

\usepackage{amssymb, latexsym}
\usepackage{color,amsmath,amssymb, amsfonts,amstext,amsthm}

\usepackage{epsfig, graphicx, graphics}
\usepackage{longtable}
\usepackage{subfigure}

\usepackage[]{hyperref}
\newcommand{\om}{\omega}

\renewcommand{\phi}{\varphi}
\newcommand{\e}{\epsilon}
\renewcommand{\a}{\alpha}

\newcommand{\R}{{\mathbb R}}

\newcommand{\eps}{\varepsilon}

\newcommand{\EX}{{\mathbb{E}}}
\newcommand{\PX}{{\mathbb{P}}}

\newtheorem{example}{Example}

\title{Quantifying model uncertainty in non-Gaussian dynamical systems
with observations on mean exit time or escape probability  }

\author{Ting Gao and Jinqiao Duan   \\
 Department of Applied Mathematics, Illinois Institute of Technology \\
  Chicago, IL 60616, USA \\
 \emph{E-mail:  tinggao0716@gmail.com,  duan@iit.edu}\\
 and \\
  Institute for Pure and Applied Mathematics, University of California\\
Los Angeles, CA 90095, USA\\
\emph{E-mail: duanjq@gmail.com }  }

\begin{document}
\date\today

\maketitle

\pagestyle{plain}

\begin{abstract}
   Complex systems are sometimes subject to non-Gaussian $\alpha-$stable L\'evy fluctuations.
   A new method is devised to estimate the uncertain parameter $\alpha$ and other system parameters, using observations on mean exit time or escape probability for the system evolution. It is based on solving an inverse problem for a deterministic, nonlocal partial differential equation via numerical optimization.
   The existing methods for estimating  parameters require observations   on system state sample paths for long time periods or   probability densities at large spatial ranges. The method proposed here,   instead,   requires   observations on mean exit time or escape probability only for an arbitrarily small spatial domain.  This new method is   beneficial to   systems for which mean exit time or escape probability is feasible to observe.

%The problem becomes a purely deterministic inverse problem.

%(ii) It also appears that this inverse problem for a nonlocal PDE with %exterior boundary condition is new (?) or at least not widely %investigated.

%  As a by-product, our work also provides a method to estimate the %$\alpha$ value (other method: estimate it by its PDF polynomial decay %rate, which requires observing probability distribution at very   %large spatial location: $x \gg 1$).

% {\bf Short Title:} Quantifying model uncertainty   \\

% {\bf Key Words:}   Stochastic dynamical systems; parameter %estimation; optimization; non-Gaussian L\'evy motion;   mean residence %time; double-well system

%{\bf Mathematics Subject Classifications (2010)}:   60H15, 60F10,
%60G17


\end{abstract}

\section{Introduction}  \label{intro}
%%%%%%%%%%%%%%%%%%%%%%%%%%%%%%%%%%%%%%%%%%%%%%%%%%%%%%%%%%%%%%%%%

Random fluctuations in   complex systems
are sometimes  non-Gaussian $\alpha-$stable L\'evy motions \cite{Woy, Shlesinger, Swinney}. We consider a   system  under such fluctuations  modeled by  a scalar  stochastic differential equation  (SDE)
\begin{equation} \label{sde}
{\rm d}X_{t}  =  f (\beta, X_{t})  {\rm d}t + \e {\rm d} L_{t}^\alpha, \;\; X_0 = x,
\end{equation}
 where $X_{t}$ is the   system state process, $f$ is a vector field (or drift), $\beta$ and $\e$ are real system parameters, and $L_{t}^\alpha$ is a scalar symmetric $\alpha-$stable L\'evy motion ($0<\alpha <2$)  defined in a probability space $(\Omega, \mathcal{F}, \mathbb{P})$. For example, the calcium
signal, as a proxy for climate state, in paleoclimatic data  is approximately described    \cite{Dit} by a model like \eqref{sde}.

 A  $\alpha-$stable L\'evy motion is a non-Gaussian process, while the well-known Brownian motion is a Gaussian process.
Non-Gaussian dynamical systems like \eqref{sde}
have attracted considerable attention recently \cite{Applebaum}, as they are appropriate models for various systems under heavy tail fluctuations \cite{Adler, taqqu}.

The process $L_t^\alpha$ has heavy tail  or power law distribution in the sense that $$
 \PX(|L_t^\alpha| > x) \sim \frac1{x^\a},
$$
for large $x$.
The $\alpha$ is called the power parameter, or stability index, or non-Gaussianity index. In fact, Brownian motion corresponds to the special case $\alpha=2$.

The $\alpha-$stable fluctuations arise  in various situations, including modeling for optimal
 foraging, human mobility and   geographical spreading of emergent infectious
disease.   GPS data are used to track the wandering
black bowed albatrosses around an island in the Southern Indian Ocean to
study their movement patterns in searching for food. It is found \cite{Humphries} that  the movement patterns
obey a power law distribution with power  parameter $\alpha \thickapprox  1.25$. One way to examine the  human mobility is to collect data by online bill trackers, which
provide   successive spatial-temporal trajectories. It is   discovered \cite{Brockmann} that the bill traveling
at certain distances within a short period of time (less than one
week) follows a  power law distribution with power parameter $\alpha \thickapprox 1.6$. Moreover, it is   noticed that the spreading patterns of human influenza, as described by    the
classic susceptibleness-infection-recovery (SIS) epidemiologic  model, is also strikingly similar to a    $\alpha-$stable L\'evy motion.

To make \eqref{sde} a predictive model, it is essential to   estimate the parameter   $\alpha$, using observations on the system evolution.   Methods for estimating other system parameters  $\beta$ and $ \e$, when $\alpha$ is known, have been considered in literature (\cite[e.g.]{HuLong1, HuLong2, JYangDuan}) and thus it is not a focus  here.
There are a couple of attempts in estimating $\alpha$.  For example, assuming the drift '$f$' insignificant (which is an inappropriate assumption in many situations), it is suggested \cite{Dit, Woy} to roughly estimate this $\alpha$ value  using data on probability density function for $X_t$.
The tail of the probability density function $p(x)$   behaves like $1/x^\alpha$ for $x \gg 1$, after ignoring   the drift '$f$'. Thus the $\log p$ vs. $\log x $ plot is a straight line with slope `$-\alpha$'. This provides an estimate $\alpha$ by   data   fitting. This method is not accurate as it assumes that the drift $f$ does not alter the tail behavior of $X_t$. Another approach to estimate $\alpha$  is suggested in \cite{Long2013} and it requires observations on lots of     system state sample   paths or  sample characteristic functions for long   time periods.

 % and thus are expensive in terms of resources, time and computation.

%The exit phenomenon, i.e., escaping from a bounded domain in state %space, is an impact of randomness on the evolution of such dynamical %systems.

\medskip

In the present paper, we devise a method to estimate $\alpha$ (and other system parameters), using observations on mean exit time or escape probability.
Recall
  the  first exit time of $X_t$ starting at $x$ (or `a particle starting at $x$')  from  a bounded  domain $D$  is defined as
\[
\tau (\omega):= \inf \{t \geq 0, X_{t}(\om, x)  \notin D \},
\]
and the mean exit time is denoted as $u(x) := \EX \tau$.
The likelihood of a
particle, starting at a point $x$, first escapes from a domain $D$ and
lands in a subset $E$ of $D^c$ (the complement of $D$) is called
escape probability  and is denoted by $P_E(x)$.

Both the mean exit time $u(x)$ and  escape probability $P_E(x)$  satisfy   deterministic, nonlocal (i.e., integral)
 differential equations with exterior Dirichlet boundary conditions.
 For the scalar SDE \eqref{sde}, these are nonlocal ordinary differential equations, while for a   SDE system, these become
 nonlocal partial differential equations.
The non-Gaussianity of the noise manifests as nonlocality at the level of the mean exit time and escape probability.

%The existing work on mean exit time gives asymptotic estimate for $u(x)$
%when the noise intensity is sufficiently small, i.e.,
%the noise term in \eqref{sde} is $\eps \; {\rm d}L_t$ with $0<\eps \ll 1$.
%See, for example, Imkeller and Pavlyukevich \cite{ImkellerP-06,
%ImkellerP-08}, and Yang and Duan \cite{YangDuan}.

When we have observations on the mean exit time $u(x)$ or  escape probability $P_E(x)$, it is thus possible to estimate $\alpha$ and other system parameters, by solving an inverse problem for the nonlocal differential equations.

It  is sometimes   too costly   to observe system state sample paths $X_t$ over very long time periods \cite{Moss}, but is more  feasible to observe (or to infer from collected data) other quantities about  system evolution, such as mean exit time and escape probability.
Mean residence time has been observed or measured in  chemical, industrial and
physiological systems \cite{Gh, Nauman}. For example,  the
mean residence time for Xe in intact and surgically
isolated muscles can be  measured \cite{Novotny}.  It is found that the mean residence time of Xe is
longer than that predicted by a single-compartment model of gas
exchange, and this leads to the  understanding that a multiple-compartment model might be more
accurate according to larger relative dispersion (the standard
deviation of residence time  divided by the mean).  Escape probability has also been observed or measured in certain physical and electronic systems \cite{Ebel, El, Forbes}.

%In the present paper, we estimate parameters $\gamma$ and $\alpha$,  %by observing the mean exit time (also called mean residence time) $u$ %in a (small) interval $D=(a, b)$.

\vskip 12pt

This paper is organized as follows. In section 2, we   formulate our method, i.e.,  an inverse problem for nonlocal differential equations to estimate parameters.
Numerical   simulation results are presented
in section 3. The paper ends with some discussions in section 4.

%%%%%%%%%%%%%%%%%%%%%%%%%%%%%%%%%%%%%%%%%%%%%%%%%
%%%%%%%%%%%%%%%%%%%%%%%%%%%%%%%%%%
%%%%%%%%%%%%%%%
%%%%%%%%%%%%%%%%%%%%%%%
%%%%%%%%%%%%%%%%%%%%%%%%%%
\section{Methods  }  \label{formulation}

A scalar symmetric $\alpha-$stable L\'evy motion $L_t^\alpha$  is characterized \cite{Applebaum, JW} by a shift coefficient which is often taken to be zero (for convenience)   and   a non-negative
measure $\nu_\a$  defined on the state space $\R^1$:
$$
\nu_\a({\rm d}x)=C_\alpha|x|^{-(1+\alpha)}\, {\rm d}x,
$$
with $\alpha\in (0, 2)$ and
$\displaystyle{C_{\alpha} =
\frac{\alpha}{2^{1-\alpha}\sqrt{\pi}}
\frac{\Gamma(\frac{1+\alpha}{2})}{\Gamma(1-\frac{\alpha}{2})}}$.
For more information see \cite{Chen, Schertzer}.
The generator for the solution process   $X_{t} $  of \eqref{sde} is
\begin{eqnarray} \label{AABB}
A  \phi &=& f(\beta, x) \phi^{'}(x)     \nonumber  \\
& &+ \e \int_{\R^1 \setminus\{0\}} [\phi(x+  y)-\phi(x) -
I_{\{|y|<1\}} \; y \phi'(x) ] \; \nu_\alpha({\rm d}y),
\end{eqnarray}
where $I_S$ is the indicator function of the set $S$, i.e.,
$$
I_S(y) =
    \begin{cases}
        1,   &\text{if $y \in S $;}\\
        0,    &\text{if $y \notin S$.}
    \end{cases}
$$

 We   consider   the mean exit time,  $u(x)$,
for an orbit starting at $x$, from a bounded open interval $D$.
By the Dynkin formula \cite{Oksendal3, Sato-99} for general Markov processes, as in \cite{Naeh, BrannanDuanErvin, BrannanDuanErvin2, Gao}, we  know that $u(x) $
satisfies the following nonlocal differential equation:
\begin{eqnarray}
 A u(x) = -1, \;\;  x\in D  \label{exit1D}\\
 u =0, \;\;  x \in D^c,     \label{exitBC}
\end{eqnarray}
where   $D^c=\R^1 \setminus D$ is the complement   of $D$.

\bigskip

%It may be    costly     to observe sample paths over very long time \cite{Moss}, but may be more  feasible to observe (or to infer from collected data) other quantities about a stochastic system, such as mean exit time \cite{Gh, Nauman} and escape probability \cite{Ebel, El, Forbes}.

Suppose that we have observed the mean exit time $u(x)$, $x\in D= (a, b)$ (a small interval). We then
solve the inverse problem for a nonlocal differential  with exterior boundary condition \eqref{exit1D}-\eqref{exitBC}, in order to estimate  $\alpha$, $\beta$ and $\e$.
%Most references consider parameter estimation in Gaussian systems %\cite{Bishwal, Ibragimov}.
  See \cite{Isakov, Belov, Kirsch} for discussions on inverse problems for partial differential equations.
This is achieved by a numerical optimization
\begin{equation}
 \min_{\alpha,\beta, \e} G(\alpha,\beta, \e), \label{optimal}
 \end{equation}
where the objective function $G= dist(u(x), u_{ob})$,  for an appropriate distance function   `dist'  between $u$ and its observation $u_{ob}$. To evaluate the objective function $G$ at initially guessed or approximated values of $(\alpha,\beta, \e)$, we need to numerically solve  \eqref{exit1D}-\eqref{exitBC} by a finite difference scheme (see Appendix).

We also consider estimation of parameters using observations on escape probability for the system \eqref{sde}. The escape probability   of a
particle, starting at a point $x$, first escapes from a bounded domain $D$ and
lands in a subset $E$ of $D^c$,    is denoted by $P_E(x)$, and it
satisfies the following nonlocal differential equation \cite{Qiao}
\begin{eqnarray}
 A\, P_E(x)&=&0, \quad x \in D, \label{eq.ep}\\
 P_E|_{x \in E}&=&1, \quad  P_E|_{x \in D^c\setminus E}=0,  \label{eq.BC}
\end{eqnarray}
where $A$ is the generator defined in \eqref{AABB}.
We  again
solve the inverse problem for a nonlocal differential  with exterior boundary condition \eqref{eq.ep}-\eqref{eq.BC}, in order to estimate  $\alpha$, $\beta$ and $\e$.
This is also achieved by a numerical optimization
\begin{equation}
 \min_{\alpha,\beta, \e} G(\alpha,\beta, \e) , \label{optimal2}
 \end{equation}
where the objective function $G = dist(P_E(x), P_{Eob})$,  for an appropriate distance    `dist'  between $P_E$ and its observation $P_{Eob}$. To evaluate the objective function $G$ at initially guessed or approximated values of $(\alpha,\beta, \e)$, we need to numerically solve  \eqref{eq.ep}-\eqref{eq.BC} by a finite difference scheme (see Appendix).

\medskip

In both settings above, the domain $D$ can be taken as small as we like (or arbitrarily small). This offers an advantage as it uses limited amount of observational resources.

In the present paper, we only consider scalar SDEs. For SDEs in higher dimensions, both mean exit time and escape probability   satisfy nonlocal partial differential  equations, and our method also applies.

%%%%%%%%%%%%%%%%%%%%%%%%%%%%%%%%%%%%%%%%%%%%%%%%%%%%%%
\section{Numerical experiments}

We now consider three examples to illustrate our method for
estimating parameters in non-Gaussian stochastic dynamical systems.
For numerical optimization, we use Matlab's
built-in function \text{fminbnd}, which is a hybrid scheme, using both
successive parabolic interpolation and golden section search to find
a   minimizer for an objective function on a fixed interval.

\begin{example}
Consider a scalar Ornstein-Uhlenbeck system
\begin{eqnarray}
dX_t = - X_t dt + dL_t^\alpha, X_0 = x.
\end{eqnarray}
In this example, $f(x)=-  x$.  Suppose that we have observed the
mean residence time $u_{ob}(x)$ for $x \in D= (-2, 2)$ and
$(-0.1,0.1)$. Let us find out estimation of $\alpha$ by solving the
inverse problem of the following nonlocal differential equation:
\begin{eqnarray}
 A u(x) = -1, \;\;  x\in D   \\
 u =0, \;\;  x \in D^c, \nonumber
\end{eqnarray}
where the generator $A$ is
\begin{eqnarray} \label{AA}
A  u &=& -x  u^{'}(x)    \nonumber \\
& &+    \int_{\R^1 \setminus\{0\}} [u(x+  y)-u(x) -
   I_{\{|y|<1\}} \, y u'(x) ] \; \nu_{\a}({\rm d}y),
\end{eqnarray}
and $D^c=\R^1 \setminus D$ is the complement set of $D$.

Using the $L^2$ norm, we define  an objective function
$$
G(\alpha)= \frac{\|u(\alpha,x) - u_{ob}(x)\|_2^2}{\|u_{ob}(x)\|_2^2},
$$
and the estimation of $\alpha \in (0, 2)$ is taken to be the minimizer, i.e.,
$$\alpha_{E} = \arg\min_{\alpha}  G(\alpha).$$
Figure \ref{OU-alpha} shows accurate  estimation of $\alpha$ on  a smaller domain  $D=(-0.1,0.1)$, as well as on a larger domain $D=(-2,2)$.

\begin{figure}[]
\begin{center}
\includegraphics*[width=6.2cm,height=5cm]{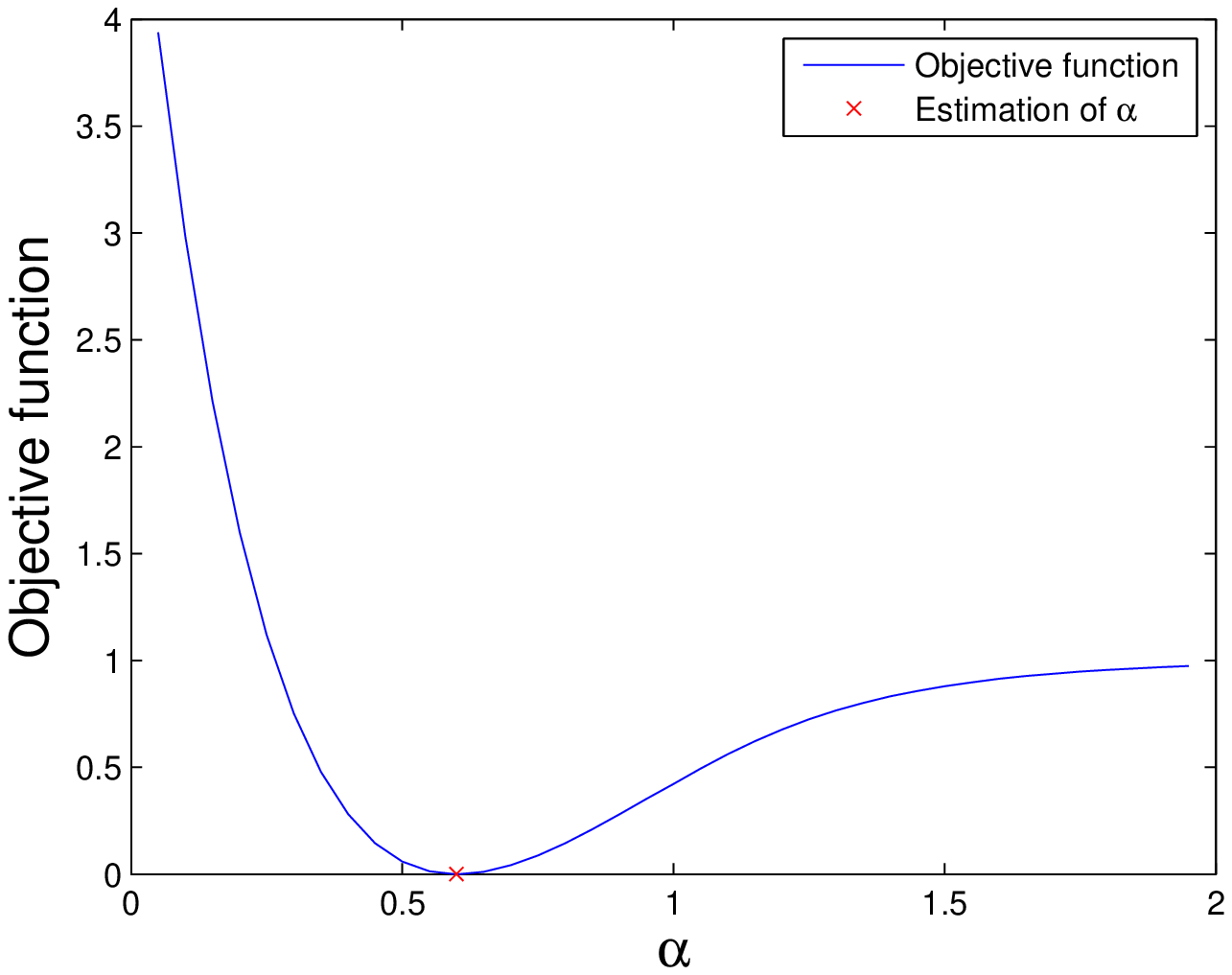}
\includegraphics*[width=6.2cm,height=5cm]{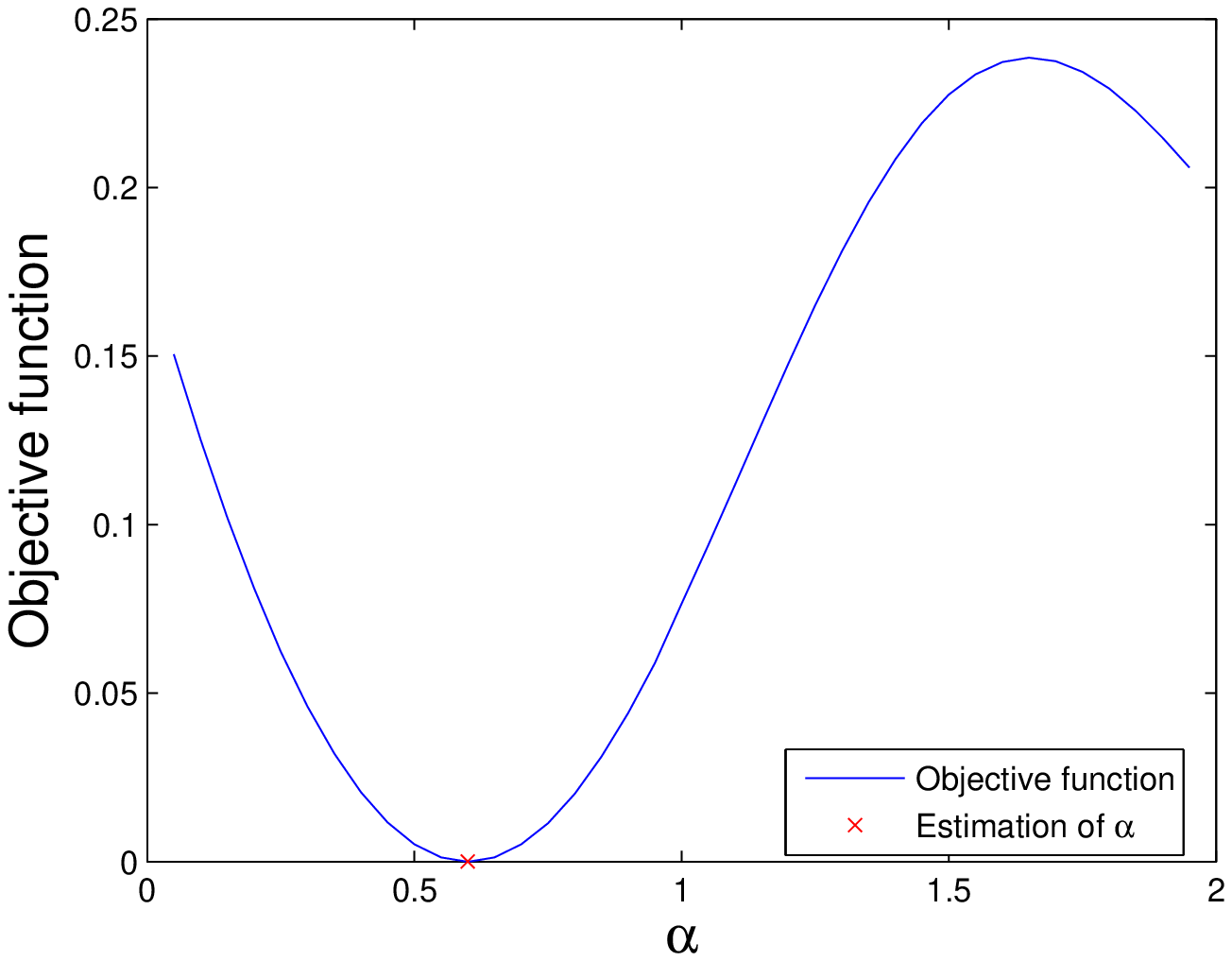}
\end{center}
\caption{Estimation of $\alpha$ on domains $D=(-0.1,0.1)$ (left) and
$D=(-2,2)$ (right) with observation on mean exit time:    True value
$\alpha=0.6$.}
\label{OU-alpha}
\end{figure}

%
%\begin{figure}[]
%\begin{center}
%\includegraphics*[width=6.2cm,height=5cm]{./figure/domain0.eps}
%\includegraphics*[width=6.2cm,height=5cm]{./figure/domain1.eps}
%\includegraphics*[width=6.2cm,height=5cm]{./figure/domain2.eps}
%\includegraphics*[width=6.2cm,height=5cm]{./figure/domain3.eps}
%\end{center}
%\caption{Estimation of different $\alpha$ along different domain
%size by mean exit time.}
%\end{figure}

%
%
%\begin{figure}[]
%\begin{center}
%\includegraphics*[width=6.2cm,height=5cm]{./figure/square.eps}
%\includegraphics*[width=6.2cm,height=5cm]{./figure/square2.eps}
%\end{center}
%\caption{Using square of norm instead of norm.}
%\end{figure}

\end{example}

\begin{example}
Consider
\begin{eqnarray}
dX_t = (X_t- X_t^3) dt +   dL_t^\alpha, X_0 = x.
\end{eqnarray}
$f(x)= x- x^3$.

We estimate $\alpha$, using observation on escape probability $P_{Eob}$.
%  We want to use escape probability $P_E(x)$
%generated from our own codes to estimate unknown parameters
%$\alpha_0$ by minimizing the distance between $P_{Eob}$ and $P_E$,
%i.e. $$\arg\min_{\alpha,\beta} dist(P_E(x,\alpha,\beta), P_{Eob}).$$
Namely, we  solve an inverse problem for the following nonlocal differential equation
\begin{eqnarray}
 A\, P_E(x)&=&0, \quad x \in D,  \\
 P_E|_{x \in E}&=&1, \quad  P_E|_{x \in D^c\setminus E}=0,
\end{eqnarray}
where $A$ is the generator defined in \eqref{AABB}.
Defining an objective function
$$
G(\alpha)= \frac{\|P_E(\alpha,x) - P_{Eob}(x)\|_2^2}{\|P_{ob}(x)\|_2^2},
$$
the estimation of $\alpha$ is $\alpha_{E} = \arg\min_{\alpha}  G(\alpha).$
Figure \ref{doublewell-alpha} shows the estimation of $\alpha=1.5$ on two different domains.

\begin{figure}[]
\begin{center}
\includegraphics*[width=6.2cm,height=5cm]{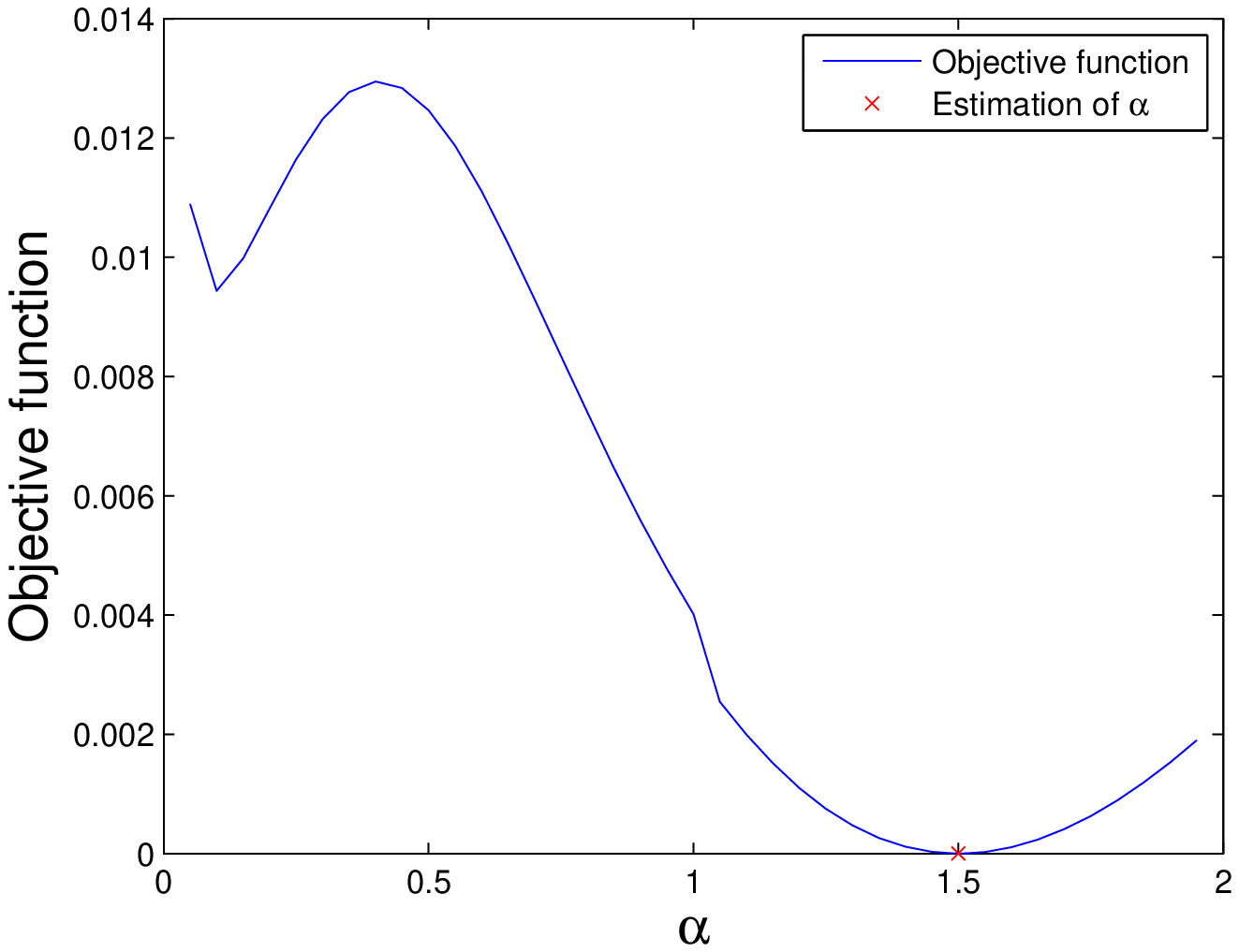}
\includegraphics*[width=6.2cm,height=5cm]{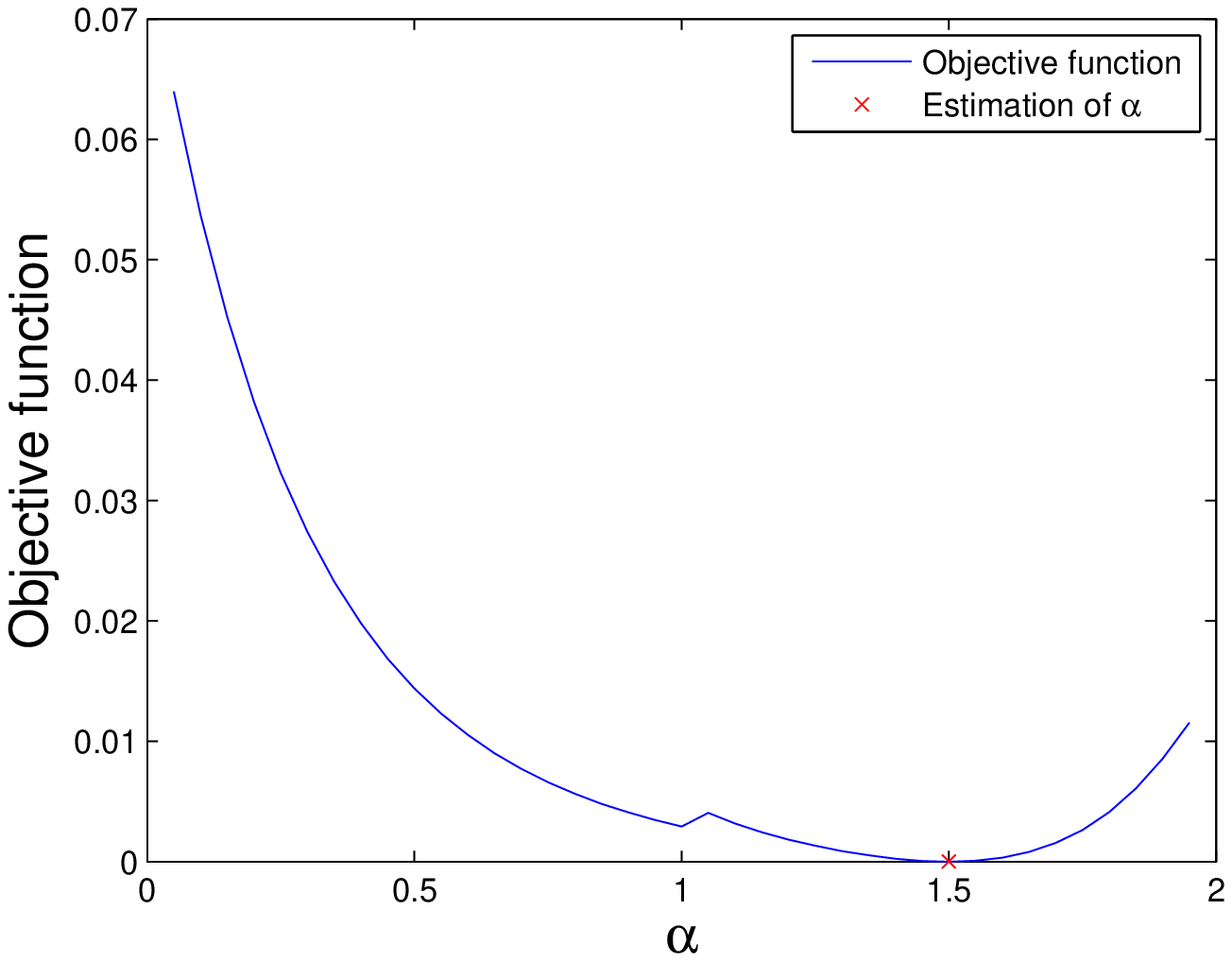}
\end{center}
\caption{Estimation of $\alpha$ on domains $D=(-0.1,0.1)$ (left) and
$D=(-2,2)$ (right) using observation on escape probability:   True value
$\alpha=1.5$.} \label{doublewell-alpha}
\end{figure}
\end{example}

\begin{example} Consider
\begin{eqnarray}
dX_t = (X_t-\beta X_t^3) dt +   dL_t^\alpha, X_0 = x.
\end{eqnarray}
$f(x)= x-\beta x^3$ where $\beta$ is a positive parameter.

In this example, we use observations of either mean exit time or escape probability
to estimate unknown parameters. Let the     observation of mean exit
time  be $u_{ob} $ and the observation of escape
probability be $P_{Eob}$.
Defining an objective function
$$
G_1(\alpha, \beta)= \frac{\|u(x,\alpha,\beta) - u_{ob}(x)\|_2^2}{\|u_{ob}(x)\|_2^2},
$$
and
$$
G_2(\alpha, \beta)= \frac{\|P_E(x,\alpha,\beta) - P_{Eob}(x)\|_2^2}{\|P_{ob}(x)\|_2^2},
$$
respectively,
we obtain estimations of parameters by minimizing these   functions separately.
Results are shown in Figure ~\ref{2dmet} for using observation of mean exit time and
Figure ~\ref{2dep} for using observation of escape probability.

%In both Fig.~\ref{2dmet} and
%Fig.~\ref{2dep}, the observation are generated when $\alpha=0.6$,
%$\beta=1.5$ and $\alpha=1.5$, $\beta=0.4$ respectively.

%Question: Is the error/residual (in 2, 1, or infinity norm) convex or not?
%Existence of global minimizer (?)

%Depending on noise intensity $\eps>0$?

%Non-convex, more minimizers?

%Interaction between nonlinearity and uncertainty.

\begin{figure}[]
\begin{center}
\includegraphics*[width=12cm,height=5cm]{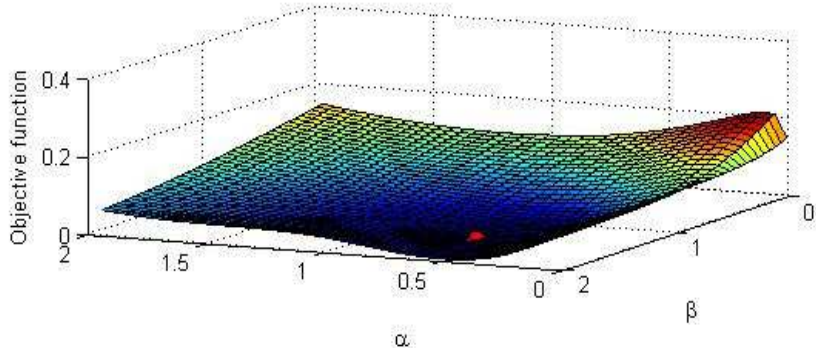}
\end{center}
\caption{Estimation of $\alpha$ and $\beta$ by observing mean exit time with
true value of $\alpha=0.6$ and true value of $\beta=1.5$. The
estimated   $\alpha$ is $0.59858$ and  the estimated
$\beta$ is $1.51$.} \label{2dmet}
\end{figure}

\begin{figure}[]
\begin{center}
\includegraphics*[width=12cm,height=6cm]{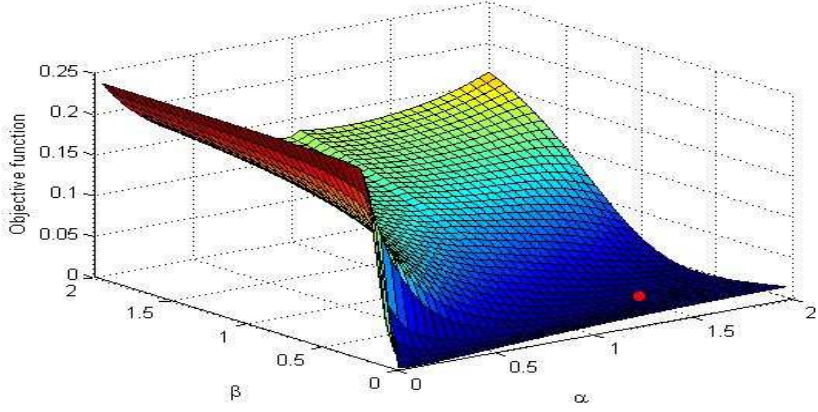}
\end{center}
\caption{Estimation of $\alpha$ and $\beta$ by observing escape probability
with true value of $\alpha=1.5$ and true value of $\beta=0.4$.  The
estimated $\alpha$ is $1.5288$  and  the estimated
$\beta$ is $0.401$.} \label{2dep}
\end{figure}

\end{example}

\section{Discussions and Conclusions}

In summary, we have devised a   method   to estimate the non-Gaussianity parameter $\alpha$,
and other system parameters, for non-Gaussian stochastic dynamical systems, using observations on either mean exit time or escape probability. It is based on solving an inverse problem for a deterministic, nonlocal partial differential equation via numerical optimization.

When the noise has a Gaussian component modeled by a Brownian motion $B_t$, the generator $A$ in nonlocal partial
differential equations \eqref{exit1D} and \eqref{eq.ep}  contains an extra Laplacian term $\Delta u$ and our method still works. Especially, if the noise has only Gaussian component,   the generator $A$ is $\Delta u$ (and the nonlocal term is absent)  and our method remains valid.

   The existing methods for estimating  the non-Gaussianity parameter $\alpha$ require observations   on system state sample paths for long time periods or   probability densities on very large spatial domain. The method proposed here,   instead,   requires   observations on either mean exit time or escape probability only for an arbitrarily small spatial domain.  This new method is especially beneficial for   systems where either mean exit time or escape probability is relatively easy to observe.

%%%%%%%%%%%%%%%%%%%%%%%%%%%%%%%%%%%%%%%%%%
%%%%%%%%%%%%%%%%%%%%%%%%%%%%%%%%%%%%%%%%%
\section*{Appendix}

In order to solve the numerical optimization problems \eqref{optimal} and \eqref{optimal2}, we need a numerical scheme to simulate the solutions of \eqref{exit1D} and \eqref{eq.ep} for given initial guesses $\alpha, \gamma, \e$, respectively. In this Appendix, we only recall a finite difference scheme \cite{Gao} for solving \eqref{exit1D}, as a similar scheme works for \eqref{eq.ep}.

Noting the principal value of
the integral  $\displaystyle{\int_\R \frac{I_{\{ |y|<\delta\}}(y)
\, y} {|y|^{1+\alpha}}\; {\rm d}y}$  always vanishes for any $\delta>0$,
we will choose the value of $\delta$ in Eq.~(\ref{exit1D})
differently according to the value of $x$.
Eq.~(\ref{exit1D}) becomes
\begin{equation} \label{exit1Dn2}
 \frac{d}{2} \; u''(x)+  f(x)\; u'(x)
  + \eps C_\a  \int_{\R \setminus\{0\}} \frac{u(x+y)-u(x)
     -  I_{\{|y|< \delta \}}(y) \; y u'(x)}{|y|^{1+\alpha}}\; {\rm d}y  = -1,
\end{equation}
for $x\in (a, b)$; and $u(x)=0$ for $x \notin  (a, b)$.

Numerical approaches for the mean exit time and escape probability in the SDEs with Brownian motions were considered in \cite{BrannanDuanErvin, BrannanDuanErvin2}, among others. In the following, we describe the numerical algorithms for the special
case of $(a,b)=(-1,1)$ for clarity of the presentation. The corresponding
schemes for the general case can be extended easily.
Because $u$ vanishes outside $(-1,1)$, Eq.~(\ref{exit1Dn2}) can be
simplified by writing $\int_{\R}=\int_{-\infty}^{-1-x} + \int_{-1-x}^{1-x}
+ \int_{1-x}^{\infty}$,
\begin{eqnarray}
  \frac{d}{2} u''(x) + f(x) u'(x)
  - \frac{\eps C_\a}{\a} \left[\frac{1}{(1+x)^\a}+\frac{1}{(1-x)^\a}\right] u(x) &
\nonumber \\
+ \eps C_\a \int_{-1-x}^{1-x} \frac{u(x+y) - u(x) -
   I_{\{|y|<\delta\}} y u'(x)}{|y|^{1+\a}}\; {\rm d}y & = -1,
\label{exit1Dn3}
\end{eqnarray}
for $x \in (-1,1)$; and $u(x)=0$ for $x \notin  (-1, 1)$.

   Noting $u$ is not smooth at the boundary points $x=-1, 1$,
in order to ensure the integrand is smooth, we rewrite Eq.~(\ref{exit1Dn3})
as
\begin{eqnarray}
  \frac{d}{2} u''(x) + f(x) u'(x)
  - \frac{\eps C_\a}{\a} \left[\frac{1}{(1+x)^\a}+\frac{1}{(1-x)^\a}\right] u(x) &
\label{exit1Dn4}\\
+ \eps C_\a \int_{-1-x}^{-1+x} \frac{u(x+y) - u(x)}{|y|^{1+\a}}\; {\rm d}y
  + \eps C_\a \int_{-1+x}^{1-x} \frac{u(x+y)-u(x)-y u'(x)}{|y|^{1+\a}}\; {\rm d}y & = -1,
\nonumber
\end{eqnarray}
for $x\geq 0$, and
\begin{eqnarray}
  \frac{d}{2} u''(x) + f(x) u'(x)
  - \frac{\eps C_\a}{\a} \left[\frac{1}{(1+x)^\a}+\frac{1}{(1-x)^\a}\right] u(x) &
\label{exit1Dn5}\\
+ \eps C_\a \int_{1+x}^{1-x} \frac{u(x+y) - u(x)}{|y|^{1+\a}}\; {\rm d}y
  + \eps C_\a \int_{-1-x}^{1+x} \frac{u(x+y)-u(x)-y u'(x)}{|y|^{1+\a}}\; {\rm d}y & = -1,
\nonumber
\end{eqnarray}
for $x < 0$. We have chosen $\delta= \text{min}\{ |-1-x|,|1-x|\}$.

   Let's divide the interval $[-2,2]$ into $4J$ sub-intervals
and define $x_j=jh$ for $-2J\leq j \leq 2J$ integer, where $h=1/J$. We denote
the numerical solution of $u$ at $x_j$ by $U_j$. Let's discretize
the integral-differential equation (\ref{exit1Dn4}) using
central difference for derivatives and ``punched-hole'' trapezoidal
rule
\begin{equation}
  \begin{split}
  & \frac{d}{2} \frac{U_{j-1} - 2U_j + U_{j+1}}{h^2}
   + f(x_j) \frac{U_{j+1} - U_{j-1}}{2h}
   -  \frac{\eps C_\a}{\a} \left[\frac{1}{(1+x_j)^\a}+\frac{1}{(1-x_j)^\a}\right] U_j \\
  & + \eps C_\a h \sum^{-J+j}_{k=-J-j}\!\!\!\!\!\!{''} \;
    {\frac{U_{j+k} - U_j}{|x_k|^{1+\a}} }
   + \eps C_\a h \sum^{J-j}_{k=-J+j,k\neq 0}\!\!\!\!\!\!\!\!\!{''} \;
    {\frac{U_{j+k} - U_j -(U_{j+1}-U_{j-1}) x_k/2h}{|x_k|^{1+\alpha}} } = -1,
  \end{split}
 \label{nm1D1a}
\end{equation}
where $j = 0,1,2, \cdots, J-1$. The modified summation symbol $\sum{''}$
means that the quantities corresponding to the two end summation
indices are multiplied by $1/2$.
\begin{equation}
  \begin{split}
  & \frac{d}{2} \frac{U_{j-1} - 2U_j + U_{j+1}}{h^2}
    + f(x_j) \frac{U_{j+1} - U_{j-1}}{2h}
   -  \frac{\eps C_\a}{\a} \left[\frac{1}{(1+x_j)^\a}+\frac{1}{(1-x_j)^\a}\right] U_j \\
  & + \eps C_\a h \sum^{J-j}_{k=J+j}\!\!\!\!{''} \;
    {\frac{U_{j+k} - U_j}{|x_k|^{1+\a}} }
   + \eps C_\a h \sum^{J+j}_{k=-J-j,k\neq 0}\!\!\!\!\!\!\!\!\!{''} \;
    {\frac{U_{j+k} - U_j -(U_{j+1}-U_{j-1}) x_k/2h}{|x_k|^{1+\alpha}} } = -1,
  \end{split}
 \label{nm1D1b}
\end{equation}
where $j = -J+1, \cdots, -2,-1$.
The boundary conditions require that
the values of $U_j$ vanish if the index $|j|\geq J$.

The truncation errors of the central difference schemes for derivatives in
(\ref{nm1D1a})  and (\ref{nm1D1b}) are of 2nd-order $O(h^2)$.
The leading-order error of the
quadrature rule is $-\zeta(\alpha-1) u''(x) h^{2-\alpha} +
O(h^2)$, where $\zeta$ is the Riemann zeta function. Thus,
the following scheme have 2nd-order accuracy
for any $0<\a<2$,
 $j = 0,1,2, \cdots, J-1$
\begin{equation}
  \begin{split}
  & C_h \frac{U_{j-1} - 2U_j + U_{j+1}}{h^2}
   + f(x_j) \frac{U_{j+1} - U_{j-1}}{2h}
   -  \frac{\eps C_\a}{\a} \left[\frac{1}{(1+x_j)^\a}+\frac{1}{(1-x_j)^\a}\right] U_j \\
  & + \eps C_\a h \sum^{-J+j}_{k=-J-j}\!\!\!\!\!\!{''} \;
    {\frac{U_{j+k} - U_j}{|x_k|^{1+\a}} }
   + \eps C_\a h \sum^{J-j}_{k=-J+j,k\neq 0}\!\!\!\!\!\!\!\!\!{''} \;
    {\frac{U_{j+k} - U_j -(U_{j+1}-U_{j-1}) x_k/2h}{|x_k|^{1+\alpha}} } = -1,
  \end{split}
 \label{nm1D2a}
\end{equation}
where $\displaystyle{C_h = \frac{d}{2} -
\eps C_\a \zeta(\alpha-1) h^{2-\a}}$.
Similarly,
for $j = -J+1, \cdots, -2,-1$,
\begin{equation}
  \begin{split}
  & C_h \frac{U_{j-1} - 2U_j + U_{j+1}}{2h^2}
    + f(x_j) \frac{U_{j+1} - U_{j-1}}{2h}
   -  \frac{\eps C_\a}{\a} \left[\frac{1}{(1+x_j)^\a}+\frac{1}{(1-x_j)^\a}\right] U_j \\
  & + \eps C_\a h \sum^{J-j}_{k=J+j}\!\!\!\!{''} \;
    {\frac{U_{j+k} - U_j}{|x_k|^{1+\a}} }
   + \eps C_\a h \sum^{J+j}_{k=-J-j,k\neq 0}\!\!\!\!\!\!\!\!\!{''} \;
    {\frac{U_{j+k} - U_j -(U_{j+1}-U_{j-1}) x_k/2h}{|x_k|^{1+\alpha}} } = -1,
  \end{split}
 \label{nm1D2b}
\end{equation}
where $j = -J+1, \cdots, -2,-1,0,1,2, \cdots, J-1$.
$U_j=0$ if $|j|\geq J$.

We solve the linear system (\ref{nm1D2a}-\ref{nm1D2b})
by direct LU factorization or
the Krylov subspace iterative method GMRES.

We find that the desingularizing term ($I_{\{|y|<\delta\}} y u'(x)$)
does not have any effect on the numerical results, regardless whether
we use LU or GMRES for solving the linear system. In this case,
we can discretize the following equation instead of \eqref{exit1Dn3}
\begin{eqnarray}
  \frac{d}{2} u''(x) + f(x) u'(x)
  - \frac{\eps C_\a}{\a} \left[\frac{1}{(1+x)^\a}+\frac{1}{(1-x)^\a}\right] u(x) &
\nonumber \\
+ \eps C_\a \int_{-1-x}^{1-x} \frac{u(x+y) - u(x)}
     {|y|^{1+\a}}\; {\rm d}y & = -1,
\label{exit1Dn6}
\end{eqnarray}
where the integral in the equation is taken as Cauchy principal value integral.
Consequently, instead of \eqref{nm1D2a} and \eqref{nm1D2b},
we have only one discretized equation for any $0<\a<2$ and
$j = -J+1, \cdots, -2,-1,0,1,2, \cdots, J-1$
\begin{equation}
  \begin{split}
  & C_h \frac{U_{j-1} - 2U_j + U_{j+1}}{h^2}
   + f(x_j) \frac{U_{j+1} - U_{j-1}}{2h} \\
  & -  \frac{\eps C_\a U_j}{\a} \left[\frac{1}{(1+x_j)^\a}+\frac{1}{(1-x_j)^\a}\right]
   + \eps C_\a h \sum^{J-j}_{k=-J-j,k\neq 0}\!\!\!\!\!\!\!\!\!{''} \;
    {\frac{U_{j+k} - U_j}{|x_k|^{1+\alpha}} } = -1.
  \end{split}
 \label{nm1D3}
\end{equation}

\bigskip

\textbf{Acknowledgement.}
 This work was partly supported by the NSF Grant 1025422. We thank Mike McCourt for help with numerical optimization.

%and the NSFC grant 11271290.

%%%%%%%%%%%%%%%%%%%%%%%%%%%%%%%%%%%%%%%%%%%%%%%%%%%%%%%%%%%%%
%%%%%%%%%%%%%%%%%%%%%%%%%%%%%%%%%%%%%%%%%%%%%%%%%%%%%%%%%%%%%%%%

\end{document}